\documentclass[12pt,a4paper]{amsart}

\usepackage{fullpage}
\usepackage[T1]{fontenc}

\usepackage{graphicx}
\newcommand{\Gr}[2]{\includegraphics[width=#2]{#1.pdf}}

\usepackage{amssymb}  
\usepackage{latexsym} 
\usepackage{colonequals}
\usepackage{comment}
\usepackage{color}
\usepackage[shortlabels]{enumitem}
\setlist[enumerate]{listparindent=0mm,parsep=0.5ex plus 0.5ex,left=0pt}
\usepackage[all]{xy}

\usepackage{charter,euler} 
\DeclareSymbolFont{bchoperators}{T1}{bch}{m}{n}
\SetSymbolFont{bchoperators}{bold}{T1}{bch}{b}{n}
\makeatletter
\renewcommand{\operator@font}{\mathgroup\symbchoperators}
\makeatother

\usepackage{titlesec} 
\titleformat{\section}{\normalfont\bfseries\filcenter}{\thesection}{1em}{}
\titleformat{\subsection}{\normalfont\bfseries}{\thesubsection}{1em}{}
\titleformat{\subsubsection}{\normalfont\bfseries}{\thesubsubsection}{1em}{}

\newcommand{\Z}{{\mathbb Z}}
\newcommand{\Q}{{\mathbb Q}}

\newcommand{\F}{{\mathbb F}}

\newcommand{\CH}{{\mathcal H}}

\newcommand{\CA}{{\mathcal A}}

\newcommand{\Gal}{\operatorname{Gal}}
\newcommand{\Aut}{\operatorname{Aut}}

\newcommand{\eps}{\varepsilon}

\newcommand{\val}{\operatorname{val}}
\newcommand{\sign}{\operatorname{sign}}
\newcommand{\Leg}[2]{\left(\frac{#1}{#2}\right)}

   {\begin{list}{}{\settowidth{\labelwidth}{#1}%
                   \setlength{\leftmargin}{\labelwidth}%
                   \addtolength{\leftmargin}{\labelsep}}}%
   {\end{list}}

\newtheorem{Theorem}{Theorem}[section]
\newtheorem{Lemma}[Theorem]{Lemma}
\newtheorem{Proposition}[Theorem]{Proposition}
\newtheorem{Corollary}[Theorem]{Corollary}
\newtheorem{Conjecture}[Theorem]{Conjecture}

\theoremstyle{definition}
\newtheorem{Definition}[Theorem]{Definition}

\newtheorem{Remark}[Theorem]{Remark}
\newtheorem{Question}[Theorem]{Question}

\numberwithin{equation}{section}

\definecolor{darkgreen}{rgb}{0,0.5,0}
\definecolor{rem}{rgb}{0.8,0,0}
\definecolor{new}{rgb}{0.3,0.1,0.9}
\definecolor{reply}{rgb}{0,0,0.8}
\definecolor{gray}{gray}{0.7}
\usepackage[
        colorlinks, citecolor=darkgreen,
        pagebackref,
        pdfauthor={Michael Stoll},
        pdftitle={When do all iterates of x² + a have maximal Galois group?},
]{hyperref}
\usepackage[alphabetic,backrefs,lite]{amsrefs} 

\begin{document}

\title{When do all iterates of $x^2 + a$ have maximal Galois group?}

\author{Michael Stoll}
\address{Mathematisches Institut,
         Universität Bayreuth,
         95440 Bayreuth, Germany}
\email{Michael.Stoll@uni-bayreuth.de}
\urladdr{\url{http://www.mathe2.uni-bayreuth.de/stoll/}}

\date{\today}


\begin{abstract}
  We provide an (in our opinion) simpler proof of the main results
  of~\cites{Stoll1992,Li2020,Li2021} on the question for which integers~$a$
  the Galois groups of all iterates of $x^2 + a$ are maximal
  (equivalently, the associated arboreal Galois representation is surjective).
  At the same time, we extend the results to rational parameters~$a$.
  The main result states that for a subset of rational numbers of
  relative density~$1/2$ in the complement of a certain open interval,
  the maximality of the Galois group of the second iterate is
  sufficient for the arboreal Galois representation to be surjective.
  We also discuss experimental evidence suggesting that a similar statement
  involving the fourth iterate may hold in general.

  The results were obtained in collaboration with Claude (Anthropic's
  AI model); they are also formalized in Lean (except for the experimental evidence).
\end{abstract}

\maketitle


\section{Introduction} 

Let $a \in \Q$ be nonzero and consider the polynomial $f = f_a = x^2 + a \in \Q[x]$.
Define its iterates by setting $f^{\circ 0} = x$ and, for $n \in \Z_{\ge 0}$,
\[ f^{\circ (n+1)} = f(f^{\circ n}) = f^{\circ n}(x^2 + a) = (f^{\circ n})^2 + a \,. \]
Assuming that no iterate has multiple roots,
we can arrange the pairs $(n, \xi)$, where $\xi$ is a root of~$f^{\circ n}$
in some fixed algebraic closure of~$\Q$, into an infinite complete binary
tree~$T$ by decreeing that $(n+1, \eta)$ is a child of~$(n, \xi)$ when
$\xi = f(\eta) = \eta^2 + a$. Then the absolute Galois group~$\Gal_{\Q}$ of~$\Q$
acts by tree automorphisms on~$T$, giving a so-called \emph{arboreal Galois representation}
\[ \rho_a \colon \Gal_\Q \to \Aut(T) \,. \]
We note that an iterate~$f^{\circ n}$ of~$f$ has multiple roots if and only if
zero is a root of~$f^{\circ m}$ for some $1 \le m \le n$. Hence some iterate
has multiple roots if and only if $0$ is periodic under~$f$. For $a \in \Q$, this
occurs exactly when $a \in \{-1, 0\}$; see Remark~\ref{R:nonzero} below.

We will discuss the following question.

\begin{Question}
  For which~$a \in \Q$ is $\rho_a$ surjective?
\end{Question}

This is equivalent to asking that the action restricted to the levels up to index~$n$,
giving the finite tree~$T_n$, be surjective for every~$n$. The automorphism
group of~$T_n$ is a Sylow $2$-subgroup of the symmetric group~$S_{2^n}$ acting
on the $2^n$~leaves of the tree. This group is isomorphic to the $n$th wreath
power of the group~$S_2$ with two elements. We will write
$\rho^n_a \colon \Gal_{\Q} \to \Aut(T_n)$ for the induced action.

This question was first asked by J.~McKay in the special case $a = 1$; see~\cite{Odoni1988}.
Odoni's earlier paper~\cite{Odoni1985} sets up some general theory.
For an overview of what is known or expected on this and related questions,
see the survey article~\cite{Jones2013} by Jones.

\begin{Definition}
  Define the sequence $(c_n)_{n \ge 0}$ (depending on~$a$) of rational numbers as follows.
  \[ c_0 = 0\,, \qquad c_1 = -a\,, \qquad c_{n+1} = f(c_n) = c_n^2 + a \quad \text{for $n \ge 1$.} \]
  This is the sequence obtained by iterating~$f$ starting with~$0$, \emph{except}
  that we have changed the sign of~$c_1$ (note that $f(-a) = f(a)$).
\end{Definition}

\begin{Remark} \label{R:nonzero}
  Note that $c_n = 0$ for some $n \ge 1$ holds if and only if $0$ is
  periodic under~$f$. For $a \in \Q$, this is the case if and only
  if $a \in \{-1, 0\}$. (This is well known; here is a quick sketch:
  if $a \notin \Z$, then the $p$-adic valuation doubles in each iteration
  step when $p$ divides the denominator, so $0$ cannot be periodic.
  So $a$ must be an integer. For $a \ge 1$ or $a \le -3$, the iteration
  sequence tends to infinity, leaving $a \in \{-2, -1, 0\}$. But for $a = -2$,
  $0$ is preperiodic ($0 \mapsto -2 \mapsto 2 \mapsto 2$) and not periodic.)
\end{Remark}

We record how the $p$-adic valuation behaves along the sequence~$(c_n)_{n \ge 1}$.

\begin{Lemma} \label{L:even_iter}
  Let $p$ be a prime and assume that $a$ is a $p$-adic integer. Then $c_n$ is
  a $p$-adic integer for all~$n$, and in~$\Z_p$, we have the congruence
  $c_{mk} \equiv c_m \bmod c_m^2$ for all $m \ge 2$ and $k \ge 1$,
  and similarly $c_k \equiv -c_1 = a \bmod c_1^2$ for all $k \ge 2$.
\end{Lemma}

\begin{proof}
  That $c_n$ is a $p$-adic integer for all~$n$ follows by an easy induction.
  Since $f = x^2 + a$ is even, so is~$f^{\circ n}$ for all $n \ge 1$,
  so $f^{\circ n}(z) \equiv f^{\circ n}(0) \bmod z^2$ for all $z \in \Z_p$.
  Also, recall that $c_m = f^{\circ m}(0)$ when $m \ge 2$ and $c_1 = -f(0)$.
  We now fix $m \ge 2$ and do induction on~$k$. The case $k = 1$ is trivial.
  For the induction step, note that
  \[ c_{m(k+1)} = f^{\circ(mk)}(c_m) \equiv f^{\circ(mk)}(0) = c_{mk} \bmod c_m^2 \,. \]
  The case $m = 1$ is similar, noting that $c_2 = f(c_1) \equiv f(0) = -c_1 \bmod c_1^2$.
\end{proof}

\begin{Lemma} \label{L:vp}
  Let $p$ be a prime and assume that $a \notin \{-1, 0\}$.
  Then exactly one of the following possibilities holds.
  \begin{enumerate}[\upshape(1)]
    \item For all $n \ge 1$, $v_p(c_n) = 0$.
    \item There are $m_p, e_p \in \Z_{\ge 1}$ such that
          \[ v_p(c_n) = \begin{cases} 0 & \text{if $m_p \nmid n$,} \\
                                      e_p & \text{if $n \ge 1$ and $m_p \mid n$.}
                        \end{cases}
          \]
    \item There is $e_p \in \Z_{<0}$ such that $v_p(c_n) = 2^{n-1} e_p$ for all $n \ge 1$.
  \end{enumerate}
\end{Lemma}

\begin{proof}
  If $v_p(c_1) = v_p(a) < 0$, then an easy induction shows that we are in the
  third case (with $e_p = v_p(a)$). If $v_p(a) \ge 0$, then we have from Lemma~\ref{L:even_iter}
  that $v_p(c_n) \ge 0$ for all~$n$. If $v_p(c_n) = 0$ for all $n \ge 1$, then we
  are in the first case. Otherwise, we are in the second case with
  $m_p = \min\{m \ge 1 : v_p(c_m) > 0\}$ and $e_p = v_p(c_{m_p})$; this follows
  from Lemma~\ref{L:even_iter} with $m = m_p$ and the observation that
  the iteration sequence modulo~$p$ and starting from~$0$ is periodic with period~$m_p$.
  See also the proof of~\cite{Stoll1992}*{Lemma~1.1(b)}.
\end{proof}

\begin{Definition}
  We say that a family $(r_i)_{i \in I}$ of rational numbers is
  \emph{$2$-independent} if for every finite subset $S \subset I$,
  $\prod_{i \in S} r_i$ is a square only when $S = \emptyset$.
  Otherwise, the family is \emph{$2$-dependent}.

  (When all $r_i$ are nonzero,
  this is equivalent to saying that the family of images of the~$r_i$
  in $\Q^\times/\Q^{\times 2}$ is linearly independent over~$\F_2$.)
\end{Definition}

Note that a family containing zero is always considered to be $2$-dependent.

The following is the equivalence ``a) $\Leftrightarrow$ b)'' in part~1) of the
Theorem in Section~1 of~\cite{Stoll1992}.

\begin{Proposition} \label{P:surj_iff}
  Fix $a \in \Q$ and $n \in \Z_{\ge 1}$. Then $\rho^n_a$ is defined and surjective
  if and only if $(c_1, \ldots, c_n)$ is $2$-independent.

  Consequently, $\rho_a$ is defined and surjective if and only if the family $(c_n)_{n \ge 1}$
  is $2$-independent.
\end{Proposition}

Note that \cite{Stoll1992} defines ``$2$-independent'' only for families
of nonzero numbers and has the general assumption that $-a$ is not a square
(which implies that all~$c_n$ are nonzero).
However, the statement remains true in the more general version stated above:
when $-a$ is a square, $\rho^1_a$ is not surjective (hence $\rho^n_a$ is not
surjective for $n \ge 1$) and $c_1 = -a$ is a square, so $(c_1, \ldots, c_n)$
is $2$-dependent. This includes the cases $a = 0$ and $a = -1$, where $\rho_a$ is not defined.

So our approach will be to show that $(c_n)_{n \ge 1}$ is $2$-independent
for $a \in \Q$ satisfying certain conditions.

The proofs in~\cites{Stoll1992,Li2020,Li2021} follow~\cite{Odoni1988}
in that they consider instead the sequence~$(b_n)$ (notation from~\cite{Stoll1992};
Odoni writes $c(n)$ and~$a(n)$ for our~$c_n$ and~$b_n$) obtained from~$(c_n)$
by multiplicative Möbius inversion. When $a \in \Z \setminus \{-1,0\}$,
then $(c_n)$ is a rigid divisibility sequence in~$\Z$ (see Lemma~\ref{L:vp}), which implies
that the~$b_n$ are pairwise coprime integers (see also \cite{Odoni1988}*{Lemma~4.4}
or~\cite{Stoll1992}*{Lemma~1.1}). The sequence~$(c_n)_{n \ge 1}$ is
$2$-independent if and only if $(b_n)_{n \ge 1}$ is, and the latter
follows if $|b_n|$ is never a square for $n \ge 2$ and $b_1$ is not a square.
Trading $(c_n)$ for~$(b_n)$ gives a simpler criterion to check --- we only
need to look at individual terms instead of arbitrary finite products ---
but at the price of having to work with more complicated objects.

One purpose of this note is to show that we can as easily (or even more easily)
obtain the results of~\cites{Stoll1992,Li2020,Li2021} on when~$\rho_a$
is surjective by working directly with~$(c_n)$.

The second purpose is to generalize the results that so far cover
certain subsets of integers~$a$ to include honest rational numbers.
The main result is that for essentially half of all rational numbers
outside the open interval~$(-2, \alpha)$, where $\alpha$ is the unique
real root of the polynomial $x^3 + 2x^2 + 2x + 2$, $\rho_a$ is surjective
as soon as (and therefore if and only if) $\rho_a^2$ is surjective.

By Proposition~\ref{P:surj_iff}, the latter means that $(c_1, c_2) = (-a, a^2 + a)$
is $2$-independent; concretely, that none of $-a$, $a^2 + a$, and $-a(a^2 + a)$
is a square. The word ``half'' above refers to six of the twelve possible
pairs of residue classes modulo~$4$ of the numerator and denominator of~$a$
(assumed to be in lowest terms with positive denominator), and ``essentially''
means that the statement holds for all such~$a$ outside~$(-2, \alpha')$,
where $\alpha' > \alpha$ is another algebraic number; for $a \in (\alpha, \alpha')$,
we can still prove it outside a density-zero subset of the half we are
considering (and for all values of~$a$ in this density-zero subset
that we checked, the claim could be verified by a finite computation).
This ``half'' can also be interpreted as a relative density:
we order the rational numbers $a = r/s$
outside the interval~$(-2, \alpha)$ by their height $\max\{|r|, s\}$;
then the proportion of the numbers of height bounded by~$X$ that are covered
by our results tends to~$1/2$ as $X$ tends to infinity.

See Section~\ref{S:rational} for precise statements.
We extend these results to two quadratic families in Section~\ref{S:quadratic};
this generalizes results from~\cite{Li2021}.
In Section~\ref{S:evidence} we add a discussion of the conjecture that
surjectivity of~$\rho_a^N$ for some fixed (sufficiently large)~$N$ should be sufficient
for surjectivity of~$\rho_a$; our new contribution is to extend the experimental
evidence to all $a \in \Q$ of height up to~$10^6$ and quadratic relations
among~$(c_1, \ldots, c_{20})$. The results support such a conjecture with $N = 4$.

A formalization in Lean (with Mathlib) of the results in this paper is available
on GitHub~\cite{GitHub}*{Part~I of the blueprint}.

\subsection*{AI declaration}

Work on this paper and its formalization was done in collaboration with Claude,
Anthropic's AI model (in various incarnations: Opus~4.8, Opus~5, Opus~5.5; Fable~5, Fable~5.1).
This proceeded in the following steps.

\begin{enumerate}[1.]
  \item When given test access to an autoformalization system by
        \href{https://www.logosresearch.ai/}{Logos Research} in early July~2026,
        I asked it to formalize~\cite{Stoll1992} (my very first published paper,
        which I had been planning to get formalized for some time).
  \item I worked quite intensely with Claude over an extended period of time on improving
        the resulting initial formalization and on including the results
        of~\cites{Li2020,Li2021}. The result is available at
        \cite{GitHub}*{Part~II of the blueprint}.
  \item This renewed my interest in the question, so
        I asked Claude to do a literature search to find out whether there
        are any further results available on the question considered here,
        in particular for rational parameters. This search came up essentially empty.
  \item I then asked Claude to see if it could use the method from~\cite{Stoll1992}
        (and its extensions by Li) to deal with more cases of integers and also
        to deal with rational parameters. Claude could extend one of Li's results
        to, in a sense, twice as many parameters. It also produced a theorem on
        rational parameters, which had a rather technical statement. At the same
        time, Claude extended the formalization to cover the new results.
  \item I then attempted to digest this by writing it up in my own words.
        While doing this, I realized that one can simplify the proofs by working
        with the sequence~$(c_n)$ directly rather than taking a detour via the
        multiplicatively Möbius-inverted sequence~$(b_n)$. I also noticed that
        one can control the sign of $c_n + c_{n+1}$ in an interval quite a bit
        to the left of~$-1$ (Claude's result required $a \le -2$ or $a > -1$).
  \item I enlisted Claude's help again. It reviewed my notes, corrected several errors, and
        contributed further ingredients that are now part of this paper: the sign pattern for
        $\alpha' < a < -1$ and criterion~\ref{T4}, the criteria \ref{Z3}, \ref{T3}
        and those at the primes dividing $r$ or~$s$, the character~$\tau$, Lemmas~\ref{L:period2}
        and~\ref{L:odd}, and the extension of the two quadratic families to rational~$\ell$.
        Prompted by me, it checked that the new approach recovers its earlier results.
        Prompted by my further questions, it found that outside $(\alpha, \alpha']$
        the conditions reduce to the $2$-independence of $c_1$ and~$c_2$ (Theorem~\ref{T:main}(1)),
        and worked out the classes in Theorem~\ref{T:main}(2)(a)--(d).
        It also formalized all results of this paper~\cite{GitHub} and made Figure~\ref{F1}.
  \item I then modified and extended the notes I had written, resulting in
        a draft of the present paper (using Claude for proof-reading).
  \item Finally, an audit of the formalization by Claude led to some further
        improvements, which I incorporated.
  \item Claude also helped with setting up a patched version of~\texttt{ratpoints}
        and running it on a number of hyperelliptic curves to obtain the experimental
        evidence in Section~\ref{S:evidence}.
\end{enumerate}

The text of this paper was written entirely by me (in some cases
using text produced by Claude as a starting point).
I take responsibility for all the claims made.


\section{Quadratic characters}

We will now show how to obtain surjectivity of~$\rho_a$ without considering
the Möbius-inverted sequence~$(b_n)$. This approach will work for rational~$a$
with certain properties; the results of~\cites{Stoll1992,Li2020,Li2021}
will follow as special cases when restricting~$a$ to be an integer.

Recall that we need to show that the family~$(c_n)_{n \ge 1}$ is
$2$-independent. To detect $2$-independence, we will use group
homomorphisms $\Q^\times \to \{\pm 1\}$.

\begin{Definition}
  We define the following homomorphisms $\Q^\times \to \{\pm 1\}$.
  When $p$ is a prime and $z \in \Q^\times$, we let $z_p$ denote an integer
  whose image in~$\Z/p\Z$ (or~$\Z/8\Z$ when $p = 2$) is the same as that of~$p^{-v_p(z)} \cdot z$.
  \begin{enumerate}[\upshape(1)]
    \item $\sign$ is the sign map;
    \item $\val_p \colon z \mapsto (-1)^{v_p(z)}$,
          where $p$ is a prime and $v_p$ is the (additive) $p$-adic valuation;
    \item $\chi_p \colon z \mapsto \bigl(\frac{z_p}{p}\bigr)$,
          where $p$ is an odd prime, i.e., $\chi_p$ is induced
          by the Legendre symbol mod~$p$ on the prime-to-$p$ part;
    \item $\chi_4 \colon z \mapsto (-1)^{(z_2 - 1)/2}$
          and $\chi_8 \colon z \mapsto (-1)^{(z_2^2 - 1)/8}$;
          these are induced by the usual Dirichlet characters modulo~$4$
          and modulo~$8$ applied to the odd part.
  \end{enumerate}
\end{Definition}

Note that the class of~$z$ modulo squares is already determined
by the values of~$\sign$ and all the~$\val_p$ on~$z$. It is, however,
helpful to use the characters that are induced by Dirichlet characters as well.

Recall from Remark~\ref{R:nonzero} that $c_n$ is nonzero for all $n \ge 1$ if and only if
$a \in \Q \setminus \{-1, 0\}$.

\begin{Definition} \label{D:phi_n}
  Let $a \in \Q \setminus \{-1, 0\}$ and $n \ge 1$.
  We say that \emph{$a$ has~$\phi_n$} if there is
  a character $\phi_n \colon \Q^\times \to \{\pm 1\}$ such that
  $\phi_n(c_n) = -1$ and $\forall m > n \colon \phi_n(c_m) = 1$.
\end{Definition}

We have the following easy consequence.

\begin{Proposition} \label{P:crit}
  Let $a \in \Q \setminus \{-1, 0\}$. Then $\rho_a$ is surjective
  if and only if $a$ has~$\phi_n$ for all $n \ge 1$.
  If $\rho_a$ is surjective, then for every sequence $(\delta_n)_{n \ge 1}$ of signs,
  there is a character $\chi \colon \Q^\times \to \{\pm 1\}$ such that
  $\chi(c_n) = \delta_n$ for all $n \ge 1$.
\end{Proposition}

\begin{proof}
  We recall that $\rho_a$ is surjective if and only if $(c_n)_{n \ge 1}$
  is $2$-independent.

  First assume that $a$ has~$\phi_n$ for all $n \ge 1$;
  we show that $(c_n)_{n \ge 1}$ is $2$-independent.
  Let $S \subset \Z_{\ge 1}$ be finite and
  nonempty and set $n \colonequals \min S$. By assumption, there
  is $\phi_n \colon \Q^\times \to \{\pm 1\}$ such that
  $\phi_n(c_n) = -1$ and $\phi_n(c_m) = 1$ for all $m > n$.
  This implies that
  \[ \phi_n\bigl(\prod_{k \in S} c_k\bigr)
       = \phi_n(c_n) \cdot \prod_{k \in S, k > n} \phi_n(c_k)
       = -1 \,;
  \]
  in particular, $\prod_{k \in S} c_k$ cannot be a square.

  Conversely, assume that $\rho_a$ is surjective, so that $(c_n)_{n \ge 1}$ is $2$-independent.
  Then the images of the~$c_n$ in~$\Q^\times/\Q^{\times 2}$ are linearly independent over~$\F_2$,
  so they can be extended to a basis. Given signs~$\delta_n$, there is therefore
  a character $\Q^\times/\Q^{\times 2} \to \{\pm 1\}$ sending the image of~$c_n$ to~$\delta_n$
  for each $n \ge 1$. This gives in particular characters that show that $a$ has~$\phi_n$
  for every $n \ge 1$.
\end{proof}

\begin{Remark}
  Alternatively, we could say that \emph{$a$ has $\phi'_n$} when there
  is a character $\phi'_n \colon \Q^\times \to \{\pm 1\}$ such that
  $\phi'_n(c_m) = 1$ for all $1 \le m < n$ and $\phi'_n(c_n) = -1$.
  Then it follows that $(c_n)_{n \ge 1}$ is $2$-independent if
  $a$ has~$\phi'_n$ for all~$n \ge 1$ (by considering $\max S$
  instead of~$\min S$). For example, $a \in \Z \setminus \{-1,0\}$ has~$\phi'_n$
  if $c_n$ has a primitive prime divisor~$p$ (i.e., a prime divisor
  not dividing any~$c_m$ for $1 \le m < n$) such that $v_p(c_n)$ is odd
  (then $\phi'_n = \val_p$ will do). It is known that for integral $a \notin \{-2,-1,0\}$,
  $c_n$ has a primitive prime divisor~$p$ whenever $n \ge 2$~\cites{Rice,DoerksenHaensch}
  (and, more generally, for rational $a = r/s \notin [-2, -1]$ except for $n = 2$
  when $r + s = 1$~\cite{Krieger}),
  but it is rather hard to directly control the parity of~$v_p(c_n)$. For example, when $a = 3$,
  we have $c_1 = -3$, $c_2 = 3 \cdot 2^2$, $c_3 = 3 \cdot 7^2$, so the
  primitive prime divisors $2$ for $n = 2$ and~$7$ for $n = 3$ occur with
  even valuation (and $c_2 c_3 = 42^2$ is a square, so $(c_n)$ is $2$-dependent
  in this case). In contrast to this, it will turn out that we are actually able to show that
  $a$ has~$\phi_n$ for all $n \ge 1$, at least for certain~$a$.
\end{Remark}

To help with constructing characters~$\phi_n$ later, it is useful to
introduce two further kinds of characters.

\begin{Definition} \label{D:phi_0 and tau}
  Let $a \in \Q \setminus \{-1, 0\}$.
  \begin{enumerate}[\upshape(1)]
    \item We say that \emph{$a$ has~$\phi_0$} if there is a character
          $\phi_0 \colon \Q^\times \to \{\pm 1\}$ such that $\phi_0(c_n) = -1$
          for all $n \ge 1$.
    \item We say that \emph{$a$ has~$\tau$} (``tail'') if there is
          a character $\tau \colon \Q^\times \to \{\pm 1\}$ such that
          $\tau(c_1) = 1$ and $\forall m \ge 2 \colon \tau(c_m) = -1$.
  \end{enumerate}
\end{Definition}

By Proposition~\ref{P:crit}, $a$ has~$\phi_0$, $\phi_1$, and~$\tau$ when $\rho_a$ is surjective.
Note that when $a$ has any two of $\phi_0$, $\phi_1$, and~$\tau$, then
$a$ also has the third (a suitable character can be obtained as the product
of the two other characters). Similarly, we have the following.

\begin{Lemma} \label{L:phi0_phin}
  Let $a \in \Q \setminus \{-1, 0\}$ and let $n \ge 1$.
  If $a$ has~$\phi_0$ and there exists
  a character $\psi \colon \Q^\times \to \{\pm 1\}$ such that
  $\psi(c_{n+m}) = -\psi(c_n)$ for all $m \ge 1$, then $a$ has~$\phi_n$.
\end{Lemma}

\begin{proof}
  We can take $\phi_n = \psi$ or $\phi_n = \phi_0  \cdot \psi$, depending
  on whether $\psi(c_n) = -1$ or~$1$.
\end{proof}

Here is a slightly more involved way of obtaining~$\phi_n$.

\begin{Lemma} \label{L:jump2}
  Let $n \ge 1$. Assume that $a \in \Q \setminus \{-1, 0\}$ has $\phi_0$,
  that some character $\psi$ satisfies $\psi(c_{m+1}) = -\psi(c_m)$ for all $m \ge n$,
  and that some character $\chi$ satisfies $\chi(c_n) = x$, $\chi(c_{n+2k}) = -x$
  for $k \ge 1$ and $\chi(c_{n+2k+1}) = y$ for $k \ge 0$, with $x, y \in \{\pm 1\}$.
  Then $a$ has $\phi_n$.
\end{Lemma}

\begin{proof}
  If $y = -x$, then $\chi$ is $x$ at $c_n$ and $-x$ at all later terms. If $y = x$,
  then the character $\chi \cdot \psi$ has the same property with $x \psi(c_n)$
  in place of~$x$. In both cases Lemma~\ref{L:phi0_phin} applies.
\end{proof}

In Section~\ref{S:phi_n}, we will use Lemma~\ref{L:phi0_phin} to show that, for the ``half''
of rational numbers mentioned in the introduction, $a$ has~$\phi_n$ for all $n \ge 2$,
provided that $a$ has~$\phi_0$. We will then consider sufficient conditions for the existence
of $\phi_0$, $\phi_1$, and~$\tau$. Combining these with Proposition~\ref{P:crit} will give
our main result.


\section{Sign information}

We record information on the sign of~$c_n$ and of $c_{n+1} + c_n$.

\begin{Lemma} \label{L:sign_cn}
  \begin{enumerate}[\upshape(1)]
    \item If $a > 0$, then $c_1 = -a < 0$ and $c_n > 0$ for all $n \ge 2$.
          In particular, $a$ has $\phi_1 = \sign$.
    \item If $-1 < a < 0$, then $c_1 = -a > 0$ and $c_n < 0$ for all $n \ge 2$.
          In particular, $a$ has $\tau = \sign$.
    \item If $a \le -2$, then $c_n \ge -a > 0$ for all $n \ge 1$.
  \end{enumerate}
\end{Lemma}

\begin{proof} \strut
  \begin{enumerate}[\upshape(1)]
    \item We have $c_{n+1} = c_n^2 + a > 0$ for all $n \ge 1$.
    \item We have $c_2 = a^2 + a = a (a + 1) < 0$. Assuming that $c_n < 0$,
          we see that $c_{n+1} = c_n^2 + a \le a^2 + a < 0$ for all $n \ge 2$;
          here we use that $c_n = c_{n-1}^2 + a \ge a$, so (using $c_n < 0$)
          $c_n^2 \le a^2$.
    \item We have $c_1 = -a > 0$. Assuming $c_n \ge -a$, we obtain
          $c_{n+1} = c_n^2 + a \ge a^2 + a \ge -a > 0$ for all $n \ge 1$.
    \qedhere
  \end{enumerate}
\end{proof}

\begin{Definition} \label{D:alpha}
  We write $\alpha$ for the unique real root of $x^3 + 2 x^2 + 2 x + 2$;
  note that $\alpha \approx -1.543689$ and that $a > \alpha$ is equivalent
  to $a^3 + 2 a^2 + 2 a + 2 > 0$.
\end{Definition}

\begin{Lemma} \label{L:sum_sign}
  Assume that $a \in (\alpha, 0)$. Then $c_{n+1} + c_n < 0$ for all $n \ge 2$.
\end{Lemma}

\begin{proof}
  First, for $-2 \le a \le 0$ we have that $c_3 \le c_2$ (since $c_3 - c_2 = a^3 (a + 2)$).
  This implies that the interval $[a, a^2 + a]$ is stable under~$f_a$.
  (Note that $f$ is convex from below, so it takes its maximum on any interval
  at one of the endpoints. Since $c_2 \ge c_3$, this maximum is~$c_2$.)
  So in particular, we have $a \le c_n \le c_2$ for all $n \ge 2$.
  The conditions $a^3 + 2 a^2 + 2 a + 2 > 0$ and $a < 0$ imply
  $c_2^2 + c_2 + a < 0$. Fix $n \ge 2$. If $c_n \ge -1/2$, then
  \[ c_{n+1} + c_n = c_n^2 + c_n + a \le c_2^2 + c_2 + a < 0 \,. \]
  Otherwise,
  \[ c_{n+1} + c_n = c_n^2 + c_n + a \le a^2 + a + a = a (a + 2) < 0 \,. \qedhere \]
\end{proof}

Note that for $a = \alpha$, we obtain $c_3 + c_2 = 0$, so the condition
in the lemma is sharp.

There is another interval where we can control the sign of~$c_n$ in a useful way.

\begin{Definition} \label{D:alpha1}
  We write $\alpha'$ for the larger real root of $x^6 + 3 x^5 + 3 x^4 + 3 x^3 + 2 x^2 + 1$.
  We have $\alpha' \approx -1.310703$ and $\alpha' < a < -1$ is equivalent to
  $c_3 < 0 < c_4$. ($c_4$ is $a (a + 1)$ times the polynomial evaluated at~$a$.)
\end{Definition}

\begin{Lemma} \label{L:sign_period2}
  Assume that $\alpha' < a < -1$. Then
  for all $k \ge 1$: $0 < c_4 \le c_{2k} \le c_2$ and $a \le c_{2k+1} \le c_3 < -1$. So
  $\bigl(\sign(c_n)\bigr)_{n \ge 1} = (1, 1, -1, 1, -1, 1, \ldots)$.
\end{Lemma}

\begin{proof}
  $c_4 = c_3^2 + a > 0$ gives $c_3^2 > -a > 1$, so $c_3 < -1$; and $c_2 = a(a+1) > 0$. On
  $B_1 = [a, c_3] \subset (-\infty, 0)$ the map $f$ is decreasing, so $f(B_1) = [c_4, c_2]$, using
  $f(a) = c_2$; on $B_2 = [c_4, c_2] \subset (0, \infty)$ it is increasing, so
  $f(B_2) = [c_5, c_3] \subseteq B_1$, as $c_5 \ge a$. Also $c_2 \in B_2$ (as $c_3 \in [a, 0)$
  gives $c_4 \le a^2 + a$) and $c_3 \in B_1$; induction on $k$.
\end{proof}

So under the assumption of the lemma, $\sign$ can serve as the character~$\psi$
in Lemma~\ref{L:jump2} when $n \ge 2$.


\section{\texorpdfstring{$\phi_n$}{phi\_n} for \texorpdfstring{$n \ge 2$}{n >= 2}}
\label{S:phi_n}

From now on, we write $a = r/s$ with $r, s \in \Z$, $s > 0$ and $\gcd(r,s) = 1$.

In this section, we will show that under suitable conditions ($r$ and~$s$
satisfy some congruence conditions modulo~$4$ and $a$ is not in a ``bad''
interval), $a$ will have~$\phi_n$ for all $n \ge 2$, assuming it also has~$\phi_0$.
The latter is needed to be able to apply Lemma~\ref{L:phi0_phin}.

The first step is to write $c_n$ as a fraction in lowest terms.
Lemma~\ref{L:vp} shows that the denominator of~$c_n$ is $s^{2^{n-1}}$ and its numerator
is divisible by~$r$. So we can write $c_n = r w_n/s^{2^{n-1}}$, where
\[ w_1 = -1 \qquad\text{and}\qquad w_{n+1} = r w_n^2 + s^{2^n-1} \]
are integers. The first few terms are
\[ w_1 = -1\,, \quad w_2 = r + s\,, \quad w_3 = r^3 + 2 r^2 s + r s^2 + s^3\,, \quad \ldots \,. \]

Since $r$ and~$s$ are coprime, $w_n$ is coprime to $r$
and~$s$, which also implies that $w_{n+1}$ and~$w_n$ are coprime.
Considered as a polynomial in $r$ and~$s$,
$w_n$ is a binary form of degree~$2^{n-1} - 1$
with integral coefficients; it is the homogenization of~$c_n/a$,
considered as a polynomial in~$a$.

We note that the numerator of $c_{n+1} + c_n$ is $r$~times
\[ v_n \colonequals r w_n^2 + w_n s^{2^{n-1}} + s^{2^n-1} = w_{n+1} + w_n s^{2^{n-1}} \,. \]
The first two terms are
\[ v_1 = r \qquad \text{and} \qquad v_2 = r^3 + 2 r^2 s + 2 r s^2 + 2 s^3 \,. \]
For $a > \alpha$ the number~$v_n$ is positive
when $n \ge 2$ (by Lemma~\ref{L:sum_sign} when $a < 0$ and obviously so when $a \ge 0$).
On the other hand, $v_n$ is negative when $a \le -2$.
Since $w_{n+1}$ and~$w_n$ are coprime, $v_n$ is coprime to~$w_n$, and since $w_{n+1}$
is coprime to~$s$, $v_n$ is also coprime to~$s$; it follows that $v_n$ is coprime
to~$w_{n+1}$ as well.

\begin{Definition}
  We set
  \[ \CA = \bigl\{a \in \Q \setminus \{-1, 0\} : \text{$a \le -2$ or $a > \alpha$}\bigr\} \,. \]
\end{Definition}

\begin{Definition}
  Let $a \in \CA$. We set
  \[ \eps = \eps_a = \begin{cases} \hfill 1 & \text{if $a > \alpha$,} \\ -1 & \text{if $a \le -2$.}
                     \end{cases}
  \]
\end{Definition}

Then the sign of~$v_n$ for $n \ge 2$ is exactly~$\eps_a$.
Note that $v_2 = s^3 (a^3 + 2 a^2 + 2 a + 2)$, which changes sign at~$\alpha$,
so we lose control over this sign when $a$ passes~$\alpha$ from above.
This is the reason why the interval~$(-2, \alpha)$ is excluded from our results.

\begin{Lemma} \label{L:phinQ}
  Let $a = r/s \in \CA$ and assume that $a$ has~$\phi_0$.
  \begin{enumerate}[\upshape(1)]
    \item \label{L:phinQ1} Fix $n \ge 2$. Then $a$ has~$\phi_n$ if $v_n$ has a prime divisor
          $p \equiv 3 \bmod 4$ (so in particular when the odd part of
          $|v_n|$ is $\equiv 3 \bmod 4$).
    \item In each of the following cases, $a$ has~$\phi_n$ for all $n \ge 2$.
          \begin{enumerate}[\upshape(C1)]
            \item \label{C1} $s$ is odd and $r + 2 s \equiv -\eps_a \bmod 4$
                  (equivalently, $r \equiv \eps_a \bmod 4$).
            \item \label{C2} $s$ is odd and $r + s \equiv -\eps_a \bmod 4$.
            \item \label{C3} $s$ is even and $r \equiv -\eps_a \bmod 4$.
          \end{enumerate}
  \end{enumerate}
\end{Lemma}

\pagebreak[1]

\begin{proof} \strut
  \begin{enumerate}[(1)]
    \item Let $p \equiv 3 \bmod 4$ be a prime divisor of $v_n = w_{n+1} + w_n s^{2^{n-1}}$.
          The considerations above then show that $p \nmid w_n, w_{n+1}, s$.
          It also follows that $p \nmid r$: $w_m \equiv s^{2^{m-1}-1} \bmod r$ for $m \ge 2$,
          so $v_n \equiv 2 s^{2^n - 1} \bmod r$.
          This implies that $c_n$ and~$c_{n+1}$ are $p$-adic units and $c_{n+1} \equiv -c_n \bmod p$.
          An easy induction then shows that $c_m \equiv c_{n+1} \bmod p$ for all $m > n$.
          Since $p \equiv 3 \bmod 4$, we have \hbox{$\chi_p(-1) = -1$}.
          Together, these observations show that $\chi_p(c_m) = -\chi_p(c_n)$ for all $m > n$.
          Now apply Lemma~\ref{L:phi0_phin}.
    \item We show that in each case, the odd part of $|v_n|$ is $\equiv 3  \bmod 4$
          for all $n \ge 2$; then part~\ref{L:phinQ1} applies.

          First assume that $s$ is odd. Then $s^{2^{m-1}} \equiv 1 \bmod 8$ for all $m \ge 2$,
          which implies that $v_n \equiv w_{n+1} + w_n \bmod 8$. If $r$ is odd, we see easily
          that $w_{2m} \equiv r + s \bmod 4$ and $w_{2m+1} \equiv s \bmod 4$ for all $m \ge 1$.
          This implies that $v_n \equiv r + 2s \bmod 4$ and hence
          that $|v_n| \equiv \eps_a (r + 2s) \bmod 4$.
          So when $r + 2 s \equiv -\eps_a \bmod 4$, it follows that $|v_n| \equiv 3 \bmod 4$.

          If $r$ is even, we find that $w_m \equiv r + s \bmod 8$ for all $m \ge 2$.
          Then $v_n \equiv 2 (r + s) \bmod 8$, hence $|v_n| \equiv \eps_a \cdot 2 (r + s) \bmod 8$.
          When $r + s \equiv -\eps_a \bmod 4$, this is $\equiv 6 \bmod 8$, so its odd
          part is $\equiv 3 \bmod 4$.

          Finally, assume that $s$ is even. Then $v_m \equiv w_{m+1} \bmod 4$ for all $m \ge 2$.
          We also have that $w_m \equiv r \bmod 4$ for all $m \ge 3$, and we conclude in the
          same way as before.
    \qedhere
  \end{enumerate}
\end{proof}

We note that for each possible sign~$\eps_a$, the three cases listed cover six
of the possible twelve primitive residue classes modulo~$4$ for pairs~$(r, s)$.
In this sense, Lemma~\ref{L:phinQ} applies to half of the rational numbers in~$\CA$.

When $(r,s)$ does not fall into one of these classes, the argument in the proof
breaks down: we obtain that $|v_n| \equiv 1 \bmod 4$ or $|v_n| \equiv 2 \bmod 8$,
which does not force the existence of a prime divisor $p \equiv 3 \bmod 4$.
It is still possible, of course, that such prime divisors exist. We will see
below that in some cases, we can show that $v_n$ (or a related number) splits
into a product whose factors can be shown to have suitable prime divisors,
but this results in a much thinner set (parameterized by a quadratic expression)
than what Lemma~\ref{L:phinQ} gives.


\section{The Klein four-group of characters}

Lemma~\ref{L:phinQ} tells us that for half of the rational $a \in \CA$,
we can show that $\rho_a$ is surjective, \emph{if} we can also produce
$\phi_0$ and~$\phi_1$. Recall that this is equivalent to $a$ having
two out of $\phi_0$, $\phi_1$, and~$\tau$: together with the trivial
character, the corresponding three characters form a Klein four-group
inside the character group of~$\Q^\times/\Q^{\times 2}$.

We will now collect sufficient conditions for each of these to be the case.
The proofs are rather straight-forward, so we keep them concise;
for details, see the Lean formalization at \cite{GitHub}*{Part~I of the blueprint}.

We note that for a quadratic character~$\psi$ on~$\Q^\times$, we have
\[ \psi(c_1) = \psi(-r s)\,, \qquad \psi(c_2) = \psi\bigl(r (r + s)\bigr) \quad\text{and}\quad
   \psi(c_1 c_2) = \psi\bigl(-s (r + s)\bigr) \,.
\]

\begin{Proposition}[Criteria for $\phi_0$] \label{prop:Z}
  $a$ has $\phi_0$ if one of the following holds.
  \begin{enumerate}[\upshape(Z1)]
    \item \label{Z1} $|r|$ is not a square.
    \item \label{Z2} A prime $p \equiv \pm 3 \bmod 8$ divides $r + 2s$.
    \item \label{Z3} A prime $p \equiv 5 \bmod 8$ divides $r^2 + s^2$ with $\chi_p\bigl((r - s) s\bigr) = -1$.
    \item \label{Z4} A prime $p \equiv 1 \bmod 4$ divides $r$ with $\chi_p(-r s) = -1$; or $8 \mid r$
          and $\chi_8(-r s) = -1$.
  \end{enumerate}
\end{Proposition}

\begin{proof}
  \ref{Z1} $v_p(c_n) = v_p(r)$ for all~$n$, so $\val_p$ is a suitable character.
  \ref{Z2} $c_n \equiv 2 \bmod p$ for all~$n$, so we can take $\chi_p$.
  \ref{Z3} Here $a^2 \equiv -1 \bmod p$, so the sequence $(c_n)_{n \ge 1}$ reduces mod~$p$
           to $(-a, a-1, -a, a-1, \ldots)$; moreover, $\chi_p(-a) = \chi_p(a) = (-1)^{(p-1)/4} = -1$
           and $\chi_p(a - 1) = \chi_p((r - s) s) = -1$. So we can take~$\chi_p$.
  \ref{Z4} We see that $\chi_p(c_1) = \chi_p(-1) \chi_p(r s)$ and
           $\chi_p(c_n) = \chi_p(r s)$ for $n \ge 2$ (and similarly for~$\chi_8$ in the second case).
           So we can take $\chi_p$, respectively,~$\chi_8$.
\end{proof}

\begin{Proposition}[Criteria for $\phi_1$] \label{prop:O}
  $a$ has $\phi_1$ if one of the following holds.
  \begin{enumerate}[\upshape(O1)]
    \item \label{O1} $a > 0$.
    \item \label{O2} $a \equiv 1 \bmod 4$, i.e., $s$ is odd and $r \equiv s \bmod 4$.
    \item \label{O3} For some prime~$p$ dividing~$s$, $-rs$ is not a square in~$\Q_p$, and
          if $p = 2$ and $v_2(s) = 2$, then moreover $\chi_4(-r s) = -1$.
    \item \label{O4} A prime $p \equiv 3 \bmod 4$ divides $r$ with $\chi_p(-r s) = -1$; or $4 \mid r$
          and $\chi_4(-r s) = -1$.
  \end{enumerate}
  In particular,  by~\ref{O3}, $a$ has~$\phi_1$ if $s$ is not a square.
\end{Proposition}

\begin{proof}
  \ref{O1} The $\sign$ character does it.
  \ref{O2} We can take $\chi_4$ or $\chi_4 \cdot \val_2$ (consideration mod~$8$). Note that
           $\val_2(c_n) = (-1)^{n-1}$.
  \ref{O3} Note that $r w_n \equiv r^{2^{n-1}} \bmod s$ for $n \ge 2$. If $v_p(s)$ is odd,
           we can take $\val_p$. If $v_p(s) \ge 2$ is even and $p$ is odd, we can take~$\chi_p$
           (the assumption on~$-rs$ implies that $\chi_p(c_1) = -1$; all $c_n$ for $n \ge 2$
           are $p$-adic squares). If $v_2(s) \ge 2$ is even, then $r w_n \equiv 1 \bmod 8$
           for $n \ge 3$ and $r w_2 \equiv 1 + r s \bmod 8$. When $v_2(s) = 2$, the additional
           assumption implies that we can take~$\chi_4$. When $v_2(s) \ge 4$, we can take
           one of $\chi_4$, $\chi_8$ and~$\chi_4 \cdot \chi_8$.
  \ref{O4} We can take $\chi_p$, resp.~$\chi_4$; compare the proof of~\ref{Z4} above.
\end{proof}

\begin{Proposition}[Criteria for $\tau$] \label{prop:T}
  $a$ has $\tau$ if one of the following holds.
  \begin{enumerate}[\upshape(T1)]
    \item \label{T1} $-1 < a < 0$.
    \item \label{T2} A prime $p \equiv 3 \bmod 4$ divides $r$ with $\chi_p(-r s) = 1$; or $4 \mid r$
          and $\chi_4(-r s) = 1$.
    \item \label{T3} A prime $p \equiv 3 \bmod 4$ divides $r + s$ with $\chi_p\bigl(-s (r+s)\bigr) = -1$; or
          $4 \mid r + s$ and $\chi_4\bigl(-s (r+s)\bigr) = -1$; or $|r + s|$ is not a square and
          divisible by $4$ or by a prime $p \equiv 3 \bmod 4$.
    \item \label{T4} $\alpha' < a < -1$ (that is, $c_3 < 0 < c_4$), and $|r + s|$ is not a square.
  \end{enumerate}
\end{Proposition}

\begin{proof}
  \ref{T1} Take the $\sign$ character.
  \ref{T2} We can take $\chi_p$ resp.~$\chi_4$; compare the proof of~\ref{Z4} above.
  \ref{T3} Let $\chi = \chi_p$, resp.~$\chi_4$.
           Here $a \equiv -1 \bmod p$, so modulo~$p$, $(c_n)_{n \ge 1}$ reduces
           to $(1, 0, -1, 0, -1, \ldots)$. We also have $(c_{2n})_p \equiv (c_2)_p \bmod p$
           for all $n \ge 1$; this follows from $c_{2k} \equiv c_2 \bmod c_2^2$.
           So $\chi$ is constant on~$(c_{2n})$. We can then take $\chi$ when $\chi(c_2) = -1$
           and $\chi \cdot \val_q$ for a prime $q$ with $v_q(r + s)$ odd otherwise.
  \ref{T4} Take $\sign \cdot \val_q$ with $v_q(r + s)$ odd.
\end{proof}

We will say that ``$a$ satisfies~(Z)'', ``$a$ satisfies~(O)'', ``$a$ satisfies~(T)''
when $a$ satisfies one of the conditions \ref{Z1}--\ref{Z4}, one of the conditions
\ref{O1}--\ref{O4}, and one of the conditions \ref{T1}--\ref{T4}, respectively.

We note a useful consequence of~\ref{Z1}, \ref{Z2} and~\ref{O2}.

\begin{Lemma} \label{L:mod4}
  If $a \equiv 1 \bmod 4$, then $a$ has~$\phi_0$ and~$\phi_1$.
\end{Lemma}

\begin{proof}
  We have that $s$ is odd and $r \equiv s \bmod 4$.
  By assumption, \ref{O2} holds. If $|r|$ is not a square, then~\ref{Z1} holds.
  Otherwise, $r = \pm t^2$ with $t$ odd,
  so $r \equiv \pm 1 \bmod 8$ and $s \equiv \pm 1 \bmod 4$ (with the same sign);
  then $r + 2s \equiv \pm 3 \bmod 8$, so~\ref{Z2} holds.
\end{proof}

One can extend these criteria to cover some further cases.
See the Lean formalization at \cite{GitHub}*{Blueprint, Section~B.1 of Part~I} for details.


\section{The main result} \label{S:rational}

As before, we write $a = r/s$ in lowest terms.

Combining the various pieces, we obtain the following result.

\begin{Theorem} \label{T:main}
  Assume that $a \in \CA$ and $(r,s)$ satisfies one of the conditions~\ref{C1}--\ref{C3}
  in Lemma~\ref{L:phinQ}. Then $\rho_a$ is surjective
  if and only if $a$ has~$\phi_0$ and~$\phi_1$. Moreover, the following hold.
  \begin{enumerate}[\upshape(1)]
    \item If $a > \alpha'$ (so $\eps_a = 1$) or $a \le -2$ (so $\eps_a = -1$),
          then the following are equivalent.
          \begin{enumerate}[\upshape(i)]
            \item $\rho_a$ is surjective;
            \item $(c_1, c_2)$ is $2$-independent;
            \item none of $-rs$, $r(r+s)$, $-s(r+s)$ is a square.
          \end{enumerate}
    \item If $\alpha < a < \alpha'$ (so $\eps_a = 1$), then $\rho_a$ is surjective
          in each of the following cases.
          \begin{enumerate}[\upshape(a)]
            \item $(r, s) \equiv (1, 1)$, $(1, 3)$ or $(3, 2) \bmod 4$.
            \item $(r, s) \equiv (2, 1) \bmod 4$ and ($r \equiv 2 \bmod 8$ or $s \equiv 5 \bmod 8$).
            \item $(r, s) \equiv (0, 3) \bmod 4$ and
                  ($|r|$ is not a square in~$\Q_2$ or $r + 2s \equiv \pm 6 \bmod 16$).
            \item $(r, s) \equiv (3, 0) \bmod 4$, $r \equiv 3 \bmod 8$ and $s \not\equiv 4 \bmod 32$.
            \item $a$ satisfies at least two of {\upshape(Z)}, {\upshape(O)} and~{\upshape(T)}.
          \end{enumerate}
  \end{enumerate}
\end{Theorem}

\begin{proof}
  If $\rho_a$ is surjective, then $a$ has~$\phi_0$ and~$\phi_1$ by Proposition~\ref{P:crit}.
  Conversely, if $a$ has~$\phi_0$ and~$\phi_1$, then $a$ has~$\phi_n$ for all $n \ge 2$
  by Lemma~\ref{L:phinQ}, so $\rho_a$ is surjective by Proposition~\ref{P:crit}.

  In~(1), the equivalence of (ii) and~(iii) and the implication ``(i) $\Rightarrow$ (ii)''
  are clear. So it suffices to show that $a$ has~$\phi_0$ and~$\phi_1$ under~(iii) in~(1)
  and in each of the cases in~(2). Since this is implied by having two out of
  $\phi_0$, $\phi_1$ and~$\tau$, case~(2e) is done; in the other cases,
  we show that condition~(2e) is also satisfied. Recall that \ref{Z1} holds if $|r|$
  is not a square, \ref{O3} if $s$ is not a square, and \ref{T3} if $|r+s|$ is not a square
  and $|r+s| \equiv 0$ or~$3 \bmod 4$. Also, when $a < -1$, so that $r < 0$ and $r + s < 0$,
  (iii) implies that at most one of~$|r|$, $s$ and~$|r + s|$ is a square, since the products
  of two of them are then~$-rs$, $r(r+s)$ and~$-s(r+s)$.

  \begin{enumerate}[\upshape(1)]
    \item \begin{enumerate}[$\bullet$]
            \item $a > 0$. \ref{O1} holds, and $r$ and $r+s$ are positive and not both squares
                  (since $r (r+s)$ is not a square by assumption).
                  In all cases except $(r,s) \equiv (1,1) \bmod 4$, this gives \ref{Z1} or~\ref{T3}.
                  In the case $(r,s) \equiv (1,1) \bmod 4$, we have $a \equiv 1 \bmod 4$,
                  and Lemma~\ref{L:mod4} applies.
            \item $-1 < a < 0$. \ref{T1} holds, and $|r|$ or $s$ is not a square,
                  giving \ref{Z1} or~\ref{O3}.
            \item $\alpha' < a < -1$. At most one of $|r|$, $s$ and $|r+s|$ is a square;
                  we obtain at least two of \ref{Z1}, \ref{O3} and~\ref{T4}.
            \item $a \le -2$. At most one of $|r|$, $s$ and $|r+s|$ is a square.
                  In all classes except $(r,s) \equiv (3,3)$, we obtain
                  at least two of \ref{Z1}, \ref{O3} and~\ref{T3}. When $(r,s) \equiv (3,3) \bmod 4$,
                  we have $a \equiv 1 \bmod 4$, and Lemma~\ref{L:mod4} applies.
          \end{enumerate}
    \item In each of the cases (a)--(d), one can show that two of (Z), (O) and~(T) hold.
          First note that $r < 0$ and $r + s < 0$.
          \begin{enumerate}[\upshape(a)]
            \item $(1,1)$: \ref{Z1} and~\ref{O2}. $(1,3)$: \ref{Z1} and~\ref{O3}. $(3,2)$: \ref{O3} and~\ref{T3}.
            \item $v_2(r) = 1$ gives~\ref{Z1}. If $s \equiv 5 \bmod 8$, then \ref{O3} holds.
                  If $r \equiv 2 \bmod 8$, then a prime $p \equiv 3 \bmod 4$ divides~$|r|$,
                  so \ref{O4} or~\ref{T2} holds.
            \item $s \equiv 3 \bmod 4$ gives~\ref{O3}. If $|r|$ is not a $2$-adic square,
                  then it is not a square in~$\Z$ and we have~\ref{Z1}.
                  If $r + 2s \equiv \pm 6 \bmod 16$, then we obtain~\ref{Z2}.
            \item $|r| \equiv 5 \bmod 8$ gives~\ref{Z1}. Here we have $v_2(s) \ge 2$, and
                  if $v_2(s) = 2$, then $s_2 \not\equiv 1 \bmod 8$ by assumption. This gives~\ref{O3},
                  with $p = 2$ unless $v_2(s)$ is even and $s_2 \equiv 5 \bmod 8$;
                  then $s$ is not a square, so \ref{O3} holds with an odd prime~$p$.
            \qedhere
          \end{enumerate}
  \end{enumerate}
\end{proof}

Since we obtain~$\phi_0$ when $|r|$ is not a square and~$\phi_1$ when $s$ is not
a square (from $\val_p$ with a prime~$p$ dividing $r$, resp.~$s$, with odd multiplicity),
we see that $\rho_a$ is surjective outside a subset of relative density zero within the
set of~$a \in \Q$ such that $\alpha < a < \alpha'$ and $a$ satisfies one of~\ref{C1}--\ref{C3}.
Even in cases where the conditions (Z), (O) and~(T) are not sufficient to conclude
that $a$ has $\phi_0$ and~$\phi_1$, we can try to find suitable characters by considering
$\chi_p$ for various~$p$ and their products. In all cases we checked, either $(c_1, c_2)$
is $2$-dependent, or else such characters can be found. There are also some further
sufficient conditions that apply in certain subfamilies. We refer to the
formalization \cite{GitHub}*{Blueprint, Appendix~B of Part~I}
for their statements and proofs. The computational observations mentioned above
suggest the following strengthening of Theorem~\ref{T:main}.

\begin{Conjecture}
  Let $a \in \CA$ satisfy one of \ref{C1}--\ref{C3}. Then $\rho_a$ is surjective
  as soon as $\rho_a^2$~is.
\end{Conjecture}

By Proposition~\ref{P:surj_iff}, $\rho_a^2$ is surjective if and only if $(c_1, c_2)$
is $2$-independent. So the conjecture holds when $a > \alpha'$ or $a \le -2$,
by Theorem~\ref{T:main}(1), and for $\alpha < a < \alpha'$, it is equivalent, by the first
statement of Theorem~\ref{T:main}, to the statement that $a$ has~$\phi_0$ and~$\phi_1$
whenever $(c_1, c_2)$ is $2$-independent.

\begin{figure}[htb]
  \Gr{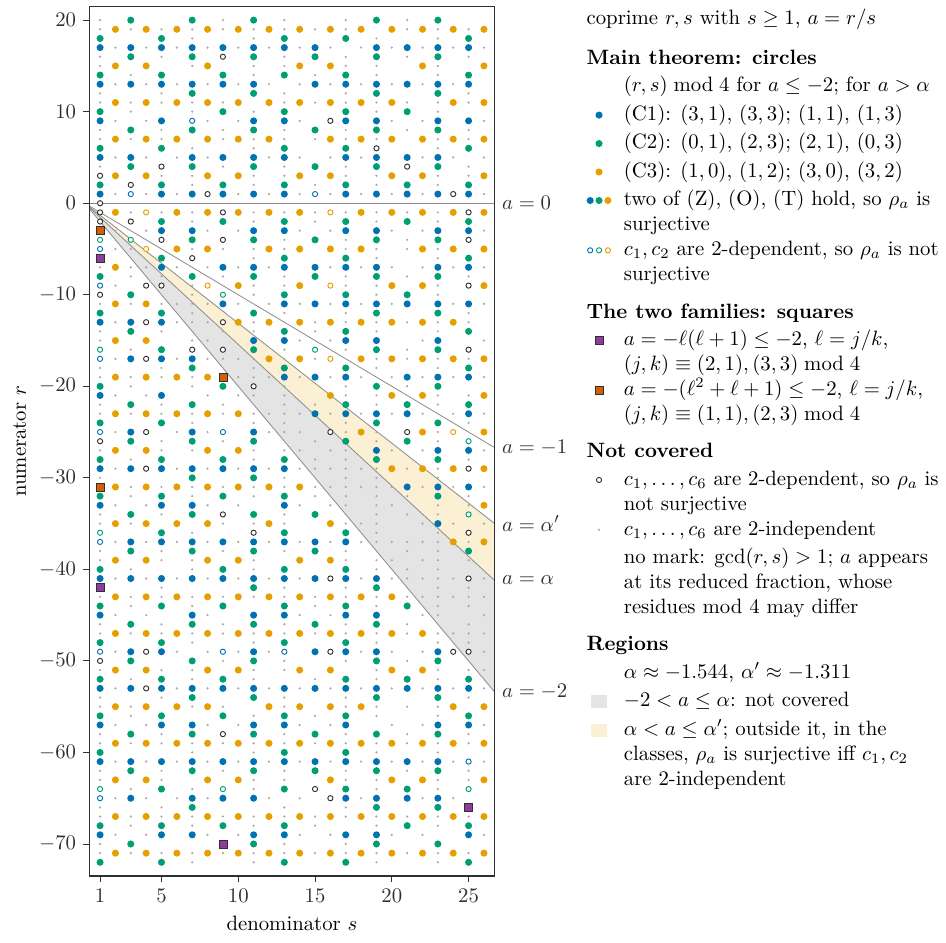}{\textwidth}
  \caption{Coverage of Theorems~\ref{T:main} and~\ref{T:quad_fam}.}
  \label{F1}
\end{figure}

Figure~\ref{F1} visualizes which rational numbers~$a$ (represented as pairs $(s,r)$ in the plane)
have $\rho_a$ surjective and are covered by Theorem~\ref{T:main} (full circles)
and Theorem~\ref{T:quad_fam} from the next section (full squares).
The hollow circles mark cases where $\rho_a$ is not surjective:
the colored ones fall under Theorem~\ref{T:main}; the black ones indicate numbers~$a$
such that the first six~$c_n$ (in fact, it is the first four in all cases in the plot)
are $2$-dependent. Gray dots represent values that are unresolved.
The shaded gray wedge corresponds to the interval~$(-2, \alpha)$; our methods
do not cover this interval.
Finally, the shaded light yellow wedge corresponds to the interval~$(\alpha, \alpha')$,
where the congruence conditions \ref{C1}--\ref{C3} are not always (known to be)
sufficient for the implication
``$(c_1, c_2)$ $2$-independent $\Rightarrow$ $\rho_a$ surjective'' to hold.

We deduce the main results of~\cite{Stoll1992} and~\cite{Li2020}.

\begin{Corollary} \label{C:integer}
  Let $a$ be an integer such that one of the following holds.
  \begin{enumerate}[\upshape(1)]
    \item $a > 0$ and $a \equiv 1$ or~$2 \bmod 4$.
    \item $a < 0$, $a \equiv 0 \bmod 4$ and $-a$ is not a square.
    \item $a < 0$, $a \equiv 3 \bmod 4$ and neither of~$-a$ and~$-a-1$ is a square.
  \end{enumerate}
  Then $\rho_a$ is surjective.
\end{Corollary}

\begin{proof}
  Since $a = -1$ is excluded by the condition ``$-a$ is not a square'', we have
  $a \le -2$ or $0 < a$; in particular, $a \in \CA$. Then, in each case, $a$
  (i.e., $(r,s) = (a,1)$) satisfies one of \ref{C1}--\ref{C3}, so part~(1)
  of Theorem~\ref{T:main} applies.

  In all cases, $c_2 = a(a+1)$ is a product of coprime nonzero integers
  differing by~$1$, hence a nonsquare. In the first case, $c_1 = -a$ and
  $c_2/c_1 = -(a + 1)$ are negative, hence nonsquares.
  In the second case, $-(a+1) \equiv 3 \bmod 4$ is a nonsquare; the other
  conditions in the last two cases are explicitly assumed.
\end{proof}


\section{Two quadratic families} \label{S:quadratic}

The main work in the proof of Theorem~\ref{T:main}, in a sense, is to show
that $a$ has~$\phi_n$ for all $n \ge 2$. The vehicle for that is the existence
of a prime $p \equiv 3 \bmod 4$ such that modulo~$p$, the sequence~$(c_n)$
becomes stationary at index~$n+1$. This follows when $|v_n| \equiv 3 \bmod 4$
or $|v_n| \equiv 6 \bmod 8$. (We could use the character induced by the
Jacobi symbol $\Leg{{\cdot}}{|v_n|}$ or $\Leg{{\cdot}}{|v_n/2|}$ to state
this in a more ``constructive'' way, i.e., without relying on a factorization.)
This is the reason for the conditions \ref{C1}--\ref{C3}: they provide the required
congruences; in all other cases the congruences fail. This is also responsible
for the fact that we only obtain a result for a set of relative density~$1/2$
in~$\Q \setminus (-2, \alpha)$.

In some cases, however, we can get a bit further. This is when we can show
that $|v_n|$ (say; there is a variant of this argument using $u_n$ as defined below),
even though congruent to~$2$ mod~$8$, systematically factors
as a product of two positive numbers one of which can be shown to be $\equiv 3 \bmod 4$
or $\equiv 6 \bmod 8$.
The main result of this section, Theorem~\ref{T:quad_fam} below, exploits this idea,
which is due to Li in~\cite{Li2021}, where essentially the special
case for integral~$a$ was obtained.

We begin by giving additional criteria that provide us with~$\phi_n$ in certain cases.
We set, for $n \ge 1$,
\[ u_n = r^2 w_n^2 - r w_n s^{2^{n-1}} + (r + s) s^{2^n - 1} \,; \]
then $c_{n+1} - c_n + 1 = u_n/s^{2^n}$ in lowest terms; in particular, $u_n$ and~$s$
are coprime. Note that
\[ c_{n+2} + c_n = (c_{n+1} + c_n) (c_{n+1} - c_n + 1) \,, \]
so modulo a prime dividing~$u_n$, we have $c_{n+2} \equiv -c_n$ (and then $c_{n+2k} \equiv -c_n$
for all $k \ge 1$ by induction).

\begin{Lemma} \label{L:period2}
  Let $n \ge 1$ and let $p \equiv 3 \bmod 4$ be a prime with $p \mid u_n$ and $p \nmid r + s$.
  Then the values of~$\chi_p$ on $(c_m)_{m\ge n}$ are $(x, y, -x, y, -x, \dots)$, with
  $x = \chi_p(c_n)$, $y = \chi_p(c_{n+1}) \in \{\pm 1\}$. Consequently, if $a$ has $\phi_0$
  and some character alternates from $n$ on, then $a$ has $\phi_n$.
\end{Lemma}

\begin{proof}
  $u_n$ is prime to $s$, so $\bar{c}_{n+2} = -\bar{c}_n$ in~$\Z/p\Z$ by the considerations
  above; it follows that
  \[ (\bar{c}_m)_{m \ge n} = (\bar{c}_n, \bar{c}_{n+1}, -\bar{c}_n, \bar{c}_{n+1}, \dots) \,. \]
  These entries are units: $\bar{c}_n = 0$ would give $p \mid r w_n$, so $p \mid r + s$,
  and $\bar{c}_{n+1} = 0$ would give $\bar{c}_n = 1$ and then $\bar{a} = -1$, so again $p \mid r + s$.
  Finally, use that $\chi_p(-1) = -1$. The last statement then follows
  from Lemma~\ref{L:jump2}, with $\chi = \chi_p$.
\end{proof}

\begin{Lemma} \label{L:odd}
  Let $a \equiv 1 \bmod 4$ with $a \le -2$. Then $a$ has $\phi_n$ for every odd $n \ge 3$.
\end{Lemma}

\begin{proof}
  Let $n \ge 3$ be odd. In $\Z/8\Z$, $\bar{c}_n = 5$ for $\bar{a} = 1$ and $\bar{c}_n = 1$ for
  $\bar{a} = 5$, and $s^{2^n} \equiv 1$, so $u_n \equiv \bar{c}_n^2 - \bar{c}_n + \bar{a} + 1 \equiv 6 \bmod 8$.
  Also $u_n > 0$: $u_n / s^{2^n} = c_n^2 - c_n + a + 1$, which is
  positive for $a \le -2$ (as $c_n \ge -a \ge 2$ gives
  $c_n(c_n - 1) + a + 1 \ge (a + 1)^2$). So a prime $p \equiv 3 \bmod 4$ divides $u_n$.
  It does not divide $r + s$: otherwise in $\Z/p\Z$
  $\bar{a} = -1$ and $(\bar{c}_m)_{m \ge 1} = (1, 0, -1, 0, -1, \dots)$, so
  $\bar{c}_n^2 - \bar{c}_n + \bar{a} + 1 = 2 \ne 0$. Lemma~\ref{L:period2}, with the alternating
  character~$\val_2$ and with~$\phi_0$ from Lemma~\ref{L:mod4}, gives~$\phi_n$.
\end{proof}

\begin{Theorem} \label{T:quad_fam}
  Let $\ell = j/k$ with coprime positive integers $j$ and~$k$.
  \begin{enumerate}[\upshape(1)]
    \item If $(j, k) \equiv (2, 1)$ or $(3, 3) \bmod 4$ and $a = -\ell(\ell + 1)$ with $a \le -2$,
          then $\rho_a$ is surjective.
    \item If $(j, k) \equiv (1, 1)$ or $(2, 3) \bmod 4$ and $a = -(\ell^2 + \ell + 1)$ with $a \le -2$,
          then $\rho_a$ is surjective.
  \end{enumerate}
\end{Theorem}

\begin{proof}
  Recall that it is sufficient to show that $a$ satisfies two of the
  conditions (Z), (O) and~(T) and that $a$ has~$\phi_n$ for all $n \ge 2$.
  \begin{enumerate}[\upshape(1)]
    \item Here $r = -j (j + k)$, $s = k^2$, $a < -2$, $\eps_a = -1$, and $v_2(r) = 1$,
          so \ref{Z1} holds. A prime $p \equiv 3 \bmod 4$ divides $r$
          (one of the positive integers $j + k$ and~$j$ is $\equiv 3 \bmod 4$), so \ref{O4} or~\ref{T2} holds.

          By Lemma~\ref{L:phinQ}~(1), $a$ has~$\phi_n$ (given that we know it has~$\phi_0$)
          when $n \ge 2$ and $|v_n|$ is divisible by a prime $p \equiv 3 \bmod 4$.
          Let $n \ge 2$, $W = -w_n$ and $V = k^{2^n - 1}$. Then
          \[ -v_n = F_1 F_2 \qquad \text{with} \qquad F_1 = (j + k) W - V \quad \text{and} \quad F_2 = j W + V \,. \]
          Since $r < 0$ and $c_n = r w_n / k^{2^n} > 0$ (Lemma~\ref{L:sign_cn}),
          we have $W > 0$, so $F_2 > 0$; then also $F_1 > 0$, as $F_1 F_2 = -v_n > 0$.
          We have $s \equiv 1 \bmod 8$, so $w_{m+1} \equiv r w_m^2 + 1$, and as $r \equiv 2 \bmod 4$
          and $w_m$ is~odd, $w_m \equiv r + 1 \bmod 8$ for $m \ge 2$; thus $W \equiv j (j + k) - 1$
          and $V \equiv k \bmod 8$. This gives $F_1 \equiv 6 \bmod 8$ when $(j, k) \equiv (2, 1) \bmod 4$
          and $F_2 \equiv 6 \bmod 8$ when $(j, k) \equiv (3, 3) \bmod 4$. In both cases,
          we obtain a prime $p \equiv 3 \bmod 4$ dividing $|v_n| = F_1 F_2$.
    \item Here $r = -(j^2 + j k + k^2)$, $s = k^2$, $r + s = -j (j + k)$, and $a \equiv 1 \bmod 4$,
          so \ref{O2} holds. $|r| \equiv 3 \bmod 4$ is not a square, so \ref{Z1} holds.
          Lemma~\ref{L:odd} gives $\phi_n$ for odd $n \ge 3$; and $\val_2$
          alternates on~$c_n$ from $n = 1$.

          Let $n \ge 2$ be even and $V = k^{2^n - 1}$. Then $u_n = G_n H_n$ with
          \[ G_n = r w_n - (j + k) V = k^{2^n}(c_n - \ell - 1) \quad \text{and} \quad
             H_n = r w_n + j V = k^{2^n}(c_n + \ell) \,,
          \]
          both positive as $c_n \ge -a = \ell^2 + \ell + 1$. By induction, we see that
          $w_n = (r + s) T$ with $T \equiv s^{2^{n-1} - 2} \bmod (r + s)$,
          so $j$ and $j + k$ divide $w_n$ (as $r + s = -j (j + k)$).
          Replacing $(j, k)$ by $(j + k, -k)$, that is, $\ell$ by $-1 - \ell$, leaves~$a$, and
          hence $r$, $s$, $w_n$ and~$T$, unchanged; it interchanges~$j$ and~$j + k$ and changes
          the sign of~$V$, so it interchanges~$G_n$ and~$H_n$, and it maps the class $(2, 3)$ to
          the class $(1, 1) \bmod 4$. The following argument uses that $j > 0$, but not the sign
          of~$k$; since $j + k > 0$, we can therefore assume that $(j, k) \equiv (1, 1) \bmod 4$.
          We set
          \[ K = H_n / j = -r (j + k) T + V > 0 \,. \]
          Then:
          \begin{itemize}
            \item $K \equiv 3 \bmod 4$.
                  In $\Z/4\Z$, $\bar{a} = 1$, so $\bar{c}_n = 2$ for even $n$, and
                  $r w_n = s^{2^{n-1}} \bar{c}_n = 2$; also $V \equiv k \bmod 4$. So
                  $j K = H_n \equiv 2 + j k$, hence $K \equiv 2 + k \equiv 3 \bmod 4$.
            \item $K$ is prime to $r + s$.
                  $K$ is odd. Let $q$ be an odd prime dividing $r + s = -j (j + k)$; then $q \nmid k$, and
                  $r \equiv -k^2$ and $r T \equiv -k^{2^n - 2} \bmod q$. So
                  $K \equiv (j + k) k^{2^n - 2} + V$, which is $2V$ if $q \mid j$ and $V$ if $q \mid j + k$.
                  The claim follows, as $q$ does not divide~$2V$.
          \end{itemize}
          So a prime $p \equiv 3 \bmod 4$ divides $K$, hence $u_n$, and $p \nmid r + s$.
          Lemma~\ref{L:period2} with $\val_2$ and $\phi_0$ gives~$\phi_n$.
    \qedhere
  \end{enumerate}
\end{proof}

Theorem~\ref{T:quad_fam} generalizes the results of~\cite{Li2021}. The first
case in the corollary below improves on~\cite{Li2021}*{Thm.~3.3} (which covers even~$m$);
the second case is~\cite{Li2021}*{Thm.~3.9}.

\begin{Corollary}
  If $m \in \Z_{\ge 0}$ and $a = -(4m + 2) (4m + 3)$ or $a = -((4m + 1) (4m + 2) + 1)$,
  then $\rho_a$ is surjective.
\end{Corollary}

\begin{proof}
  Apply Theorem~\ref{T:quad_fam} with $k = 1$.
\end{proof}


\section{Computational evidence} \label{S:evidence}

In this section, we report on computational evidence supporting the following conjecture.

\begin{Conjecture} \label{C:finite-infinite}
  There is a natural number~$N$ such that for all $a \in \Q$, if $\rho_a^N$ is surjective,
  then so is~$\rho_a$.
\end{Conjecture}

Forms of this conjecture are already present in the literature.
Hindes~\cite{Hindes2015} considers the set~$S^{(n)}$ of $a \in \Q$ such that
$\rho_a^{n-1}$ is surjective, but $\rho_a^n$ is not. At the end of Section~3, he raises
the question whether $S^{(n)} = \emptyset$ for all sufficiently large~$n$, which is
equivalent to Conjecture~\ref{C:finite-infinite}, and asks ``Is it true for all $n \ge 5$?'',
which is Conjecture~\ref{C:finite-infinite} with $N = 4$.

Our main result (Theorem~\ref{T:main}) shows that Conjecture~\ref{C:finite-infinite}
holds for half of the rational numbers outside the interval $(-2, \alpha')$
with $N = 2$ (and within the interval $(\alpha, \alpha')$ for half of the numbers
outside at most a set of density zero).
This rests on the fact that such $a$ have $\phi_n$ for all $n \ge 2$ (assuming
they have~$\phi_0$), which follows from the fact that $a \in \CA$ and satisfies
one of the congruence conditions~\ref{C1}--\ref{C3}. Our method cannot reach
much beyond this (except for the thin quadratic families discussed in Section~\ref{S:quadratic}).

So, to investigate whether Conjecture~\ref{C:finite-infinite} has a chance
of being true, we tested all rational numbers~$a$ of height up to~$10^6$
(i.e., $a = r/s$ with $|r| \le 10^6$ and $1 \le s \le 10^6$) for quadratic relations
among $c_1, \ldots, c_{20}$, i.e., we determined the non-empty subsets~$T$
of~$\{1, \ldots, 20\}$ with $\max T \ge 4$ such that $\prod_{n \in T} c_n$ is a square in~$\Q$.
This gives the following result.

\begin{Proposition} \label{P:evidence}
  Let $a \in \Q$ be of height at most~$10^6$. Assume that $\rho_a^3$ is surjective and
  $a \notin \{2/3, -6/7\}$. Then $\rho_a^{20}$ is surjective.
  For $a \in \{2/3, -6/7\}$, $\rho_a^3$ is surjective, but $\rho_a^4$ is not.
\end{Proposition}

\begin{proof}
  Recall that by Proposition~\ref{P:surj_iff}, $\rho_a^n$ is surjective if and only
  if $(c_1, \ldots, c_n)$ is $2$-independent.

  Let $C_0 = 0 \in \Z[x]$ and $C_{n+1} = C_n^2 + x$ for $n \ge 0$ (so that
  for any given~$a$, $c_1 = -C_1(a)$ and $c_n = C_n(a)$ for $n \ge 2$). There
  are monic polynomials $B_n \in \Z[x]$ for $n \ge 1$ such that
  $C_n = \prod_{d \mid n} B_d$ for all $n \ge 1$, and the pairwise resultants
  of the~$B_n$ are all $\pm 1$; see~\cite{BEK2020}*{Lemma~9.5}, which the authors
  attribute to Bjorn Poonen. (The polynomials~$C_n$ are known as Gleason polynomials;
  the complex roots of the~$B_n$ are the centers of period exactly~$n$ in the Mandelbrot set.)

  By Remark~\ref{R:nonzero}, the only rational roots of any of the~$C_n$ and
  therefore of any of the~$B_n$ are $0$ and~$-1$, for which already~$\rho_a^1$
  is not surjective, so that the claim holds for them. We also know that
  $\rho_{-2}^2$ is not surjective ($c_1 c_2 = 2^2$ in this case).
  We can therefore assume from now on that $a \notin \{-2, -1, 0\}$.
  So let $\CH$ denote the set of all rational numbers~$a$
  of height at most~$10^6$ except $a = -2$, $a = -1$, and $a = 0$.
  Since no factor in a product $\prod_{n \in T} c_n$ vanishes,
  the quadratic relations among the~$c_n$ correspond one-to-one to quadratic relations
  between the numbers $B_n(a)$, except for a sign coming from the fact
  that $c_1 = -C_1(a)$, which is therefore present if and only if the relation of the~$c_n$
  involves~$c_1$; this is equivalent to the degree of the product of~$B_n$ obtained
  by canceling square factors from the corresponding product of the~$C_n$ being odd
  (note that $\deg C_n = 2^{n-1}$ is even for $n \ge 2$). It is not hard to show
  by induction that $B_n(-1) = B_n(0) = 1$ and $B_n(-2) = (-1)^{\deg B_n}$ for all $n \ge 3$.
  This implies that even-degree products of the~$B_n$ with $n \ge 3$ evaluated at~$a$
  are $3$-adic squares for all $a \in \Q$; in particular, the negative of such
  a product can never be a square in this case.

  The fact that the pairwise resultants of the~$B_n$ are units implies that
  in any product of values~$B_n(a)$ (for $a \in \Q$ and pairwise distinct indices~$n$)
  that is a square up to a sign,
  each value of a subproduct such that the corresponding product of the~$B_n$
  has even degree must be a square up to a sign. This follows from considering
  the even-degree homogenization of the product polynomial evaluated at
  the numerator and denominator of~$a$: the unit resultants imply coprimality
  of the factors coming from the various~$B_n$ (and, as the $B_n$ are monic,
  with the denominator as well). When the degree of the product
  polynomial is odd, this argument also shows that the value of each odd-degree
  factor must be a square up to sign. The $3$-adic consideration above
  then lets us exclude the negative sign in each case where we have an
  even-degree product not involving $B_1$ or~$B_2$. Also note that the degree of~$B_n$
  is odd if and only if $n$ is squarefree. From these considerations we deduce that the claim
  follows if we can verify the following.
  \begin{enumerate}[(1)]
    \item For $8 \le n \le 20$ not squarefree and $a \in \CH$, $B_n(a)$ is not a square.
    \item For $a \in \CH$, $B_4(a)$ is a square if and only if $a = 2/3$ or $a = -6/7$.
    \item For $m, n$ squarefree with $3 \le m < n \le 19$ and $a \in \CH$, $B_m(a) B_n(a)$
          is not a square.
    \item For $5 \le n \le 19$ squarefree and $a \in \CH$, both $B_1(a) B_n(a)$
          and~$B_2(a) B_n(a)$ are not of the form~$\pm \square$.
    \item For $5 \le n \le 19$ squarefree and $a \in \CH$, $B_n(a)$ is not of the form~$\pm\square$.
  \end{enumerate}
  (Note that these imply that for $a \in \CH$, there cannot be any
  relation $\prod_{n \in T} c_n = \square$ with $5 \le \max T \le 20$, and there is
  such a relation with $\max T = 4$ only for $a = 2/3$ and $a = -6/7$.)

  Each of these claims is equivalent to saying that on one or several hyperelliptic curves,
  all finite rational points whose $x$-coordinate has height at most~$10^6$
  have $x$-coordinate in $\{-2, -1, 0\}$ (or in~$\{-2, -6/7, -1, 0, 2/3\}$ in the case $n = 4$).
  The \texttt{ratpoints} program~\cite{ratpoints} (written and maintained by the author)
  can search for rational points on hyperelliptic curves very efficiently.
  With some patches to enable it to deal with polynomials of large degrees
  ($\deg B_{20} = 523\,770$) and with large coefficients, we ran it on all these curves
  with a height bound of~$10^6$. This took about three hours; the most
  time-intensive part was to generate the more than~170\,GB of coefficient
  data for the polynomials. The results confirm that the only rational points
  up to that height on these curves are the expected ones.

  For $a = 2/3$ and~$-6/7$, we easily check that $(c_1, c_2, c_3)$ is $2$-independent,
  whereas $c_2 c_4$ is a square in both cases (this is equivalent to $B_4(a)$ being
  a square). This implies the claims in the last sentence of the statement
  of Proposition~\ref{P:evidence}.
\end{proof}

\begin{Remark}
  Consider $\prod_{n \in T} B_n = (-1)^t y^2$, where $t = \sum_{n \in T} \deg B_n$ is
  the degree of the left-hand side. (Note that here $T$ denotes a set of indices
  in a product of the~$B_n$, rather than in a product of the~$c_n$ as in the proof above.
  The sign $(-1)^t$ is negative if and only if the corresponding product of~$C_n$'s
  involves~$C_1$; recall that $C_1(a) = a = -c_1$.)
  For $T = \{1\}$, $\{2\}$, $\{1, 2\}$, the
  resulting equation gives a curve of genus~$0$ whose rational points are easily
  parameterized. As $\deg B_3 = 3$, we obtain curves of genus~$1$ for $T = \{3\}$,
  $\{1, 3\}$, and~$\{2, 3\}$. The first and the third of these are elliptic curves of rank~$1$ and
  therefore lead to infinitely many rational points, whereas the second one has rank~$0$
  and rational points only at infinity and with $x = 0$.
  For $T = \{1, 2, 3\}$ we obtain a genus~$2$ curve whose Jacobian has a finite group
  of rational points, and we easily obtain that its rational points are at infinity
  or have $x = 0$ or~$-1$. For these results, see also~\cite{Hindes2013}.
  For $T = \{4\}$ we obtain another genus~$2$ curve.
  Its Jacobian has Mordell--Weil rank~$3$, so Chabauty's method does not apply,
  and the Galois group of~$B_4$ is maximal (i.e.,~$S_6$), so alternative ways of determining the
  rational points are not very promising; see also~\cite{Hindes2015}, where also
  a computation is described that shows that there are no further rational points
  up to height~$10^{100}$ by running a Mordell--Weil sieve. All other curves have genus at least~$3$;
  Hindes in \cite{Hindes2018}*{Lemmas~4.2, 4.3} uses the Chabauty--Coleman method
  to show that there are no further rational points on $y^2 = \pm B_5(x)$ (these
  are curves of genus~$7$ whose Jacobians have Mordell--Weil ranks $5$ and~$3$).
  In the preprint~\cite{Stoll2019}, we show that the curve $y^2 = C_7(x) = B_1(x) B_7(x)$
  of genus~$31$ has only the points at infinity and~$(0, 0)$ as rational points
  (assuming GRH for the computation of the class group of the number field of
  degree~$63$ obtained by adjoining a root of~$B_7$ to~$\Q$). We also show
  that the curves $y^2 = \pm C_5(x) = \pm B_1(x) B_5(x)$ of genus~$7$ have
  no unexpected rational points.
\end{Remark}

\begin{Remark}
  Restricting the parameter~$a$ to integer values, one can show a bit more.
  Hindes in~\cite{Hindes2018}*{Thm.~1.3} shows that in this case, the surjectivity
  of~$\rho_a^3$ implies that of~$\rho_a^5$ and conjectures that this also
  implies that~$\rho_a$ is surjective (Conj.~1.5 in \emph{loc.~cit.}). This is
  Conjecture~\ref{C:finite-infinite} for integral~$a$ with $N = 3$.
  In~\cite{Hindes2016}, he shows that the integral version of the conjecture
  follows from the Vojta or the Hall--Lang conjecture.
\end{Remark}

\begin{Remark}
  One can consider a more general version of the question we are studying in this paper,
  by replacing~$0$ with another rational base point~$b$ and looking at the binary
  tree of iterated preimages of~$b$ under~$f_a$. This leads to arboreal Galois representations
  $\rho_{a,b}^n$ and~$\rho_{a,b}$. Equivalently, one considers the shifted polynomial
  $(x + b)^2 + a - b$. In this setting, there is an example
  found by Jones~\cite{Jones2012}*{Eq.~(1)}, where $\rho_{a,b}^8$ is surjective,
  but $\rho_{a,b}^9$ is not. (This example is $x^2 + 1$ with $b = -88\,255\,775\,491\,812\,351\,975\,604$.)
\end{Remark}

If we assume that no further (and thus fairly large) points exist on these curves,
then the proof of Proposition~\ref{P:evidence} implies that if Conjecture~\ref{C:finite-infinite}
holds for some $N \le 20$, then it also holds for $N = 4$. This would suggest
that the conjecture will likely indeed hold with $N = 4$. (The conjecture is wrong
with $N < 4$, as the examples $a = 2/3$ and~$-6/7$ show.)


\end{document}